\documentclass[a4paper, 11pt]{amsart}
\usepackage[latin1]{inputenc}
\usepackage[ T1]{fontenc}
\usepackage[english]{babel}
\usepackage{amssymb}
\usepackage{amsmath}
\usepackage{amsthm}
\usepackage{amscd}
\usepackage{amsfonts}
\usepackage{stmaryrd}
\usepackage{pb-diagram}
\usepackage{epic,eepic,epsfig}
\usepackage{a4wide}
\usepackage{nextpage}
\usepackage{fancyhdr}
\usepackage{enumerate}
\usepackage{hyperref}
\usepackage{float}
\usepackage{tikz}

\usepackage{enumitem} 
\setlist[enumerate]{label=\textnormal{(\arabic*)}}

\newtheorem{prop}{Proposition}[section]

\newtheorem{lemma}[prop]{Lemma}

\newtheorem{theorem}[prop]{Theorem}

\newtheorem{proposition}[prop]{Proposition}

\theoremstyle{definition}
\newtheorem{rem}[prop]{Remark}

\numberwithin{equation}{section}

\newcommand{\PP}{{\mathbb{P}}}
\newcommand{\QQ}{{\mathbb{Q}}}
\newcommand{\RR}{{\mathbb{R}}}
\newcommand{\ZZ}{{\mathbb{Z}}}
\newcommand{\CC}{{\mathbb{C}}}

\newcommand{\Spec}{\operatorname{Spec}}

\newcommand{\OO}{{\mathcal{O}}}

\newcommand{\ndot}{\raisebox{.4ex}{.}}

\newcommand{\Proof}{{\sl Proof.}\quad}
\newcommand{\QED}{{\unskip\nobreak\hfil\penalty50\quad\null\nobreak\hfil
{$\Box$}\parfillskip0pt\finalhyphendemerits0\par\medskip}}

\usepackage{colonequals}

\title{Absolute minimum of a height associated to a semipositive adelic line bundle}
\author{Shu Kawaguchi}
\address{Department of Mathematics, Graduate School of Science, Kyoto University, Kyoto 606-8502, Japan.}
\email{kawaguch@math.kyoto-u.ac.jp}
\author{Fabien Pazuki}
\address{Department of Mathematical Sciences, University of Copenhagen,
Universitetsparken 5, 
2100 Copenhagen \O, Denmark.}
\email{fpazuki@math.ku.dk}

\date{July 21, 2026}

\begin{document}

\begin{abstract}
There exists a semipositive adelic line bundle on a projective variety over $\QQ$ whose absolute minimum does not attain the minimum by algebraic points. The construction is based on arithmetic properties of a family of dynamical systems over $\mathbb{P}^1$.
\end{abstract}

\maketitle

\smallskip

    
\smallskip

\noindent{\bf 2020 Math. Subj. Classification:} 11G50, 14G40, 37P30.



\section{Introduction}
Let $K$ be a number field and denote by $M_K$ the set of 
places of $K$. Let $X$ be a normal and geometrically integral projective variety over $K$ and let $\overline{D} = (D, (g_v)_{v \in M_K})$ be a semipositive adelic $\RR$-Cartier divisor on $X$ such that $D$ is ample. For a subvariety $Y \subseteq X$, let $\widehat{h}_{\overline{D}}(Y)$ denote the normalized height of $Y$ associated to $\overline{D}$. The absolute minimum of $\overline{D}$ is defined as~$\zeta_\mathrm{abs}(\overline{D}) \colonequals \inf_{x \in X(\overline{K})} \widehat{h}_{\overline{D}}(x)$. For the notion of adelic $\RR$-Cartier divisors, we refer to Moriwaki \cite{Mo16}, which generalizes the notion of adelic line bundles due to Zhang \cite{Zh95b}.

In \cite[Theorem~1.2]{Ba23}, Balla\"{y} showed that the inequality $\zeta_\mathrm{abs}(\overline{D}) \leqslant \widehat{h}_{\overline{D}}(Y)$ holds for any subvariety $Y \subseteq X$ and that equality is attained for some subvariety $Y_0 \subseteq X$. In Remark~5.4 of {\it loc.\! cit.}, an anonymous referee raised the question of whether $Y_0$ can be chosen to be zero-dimensional, or equivalently, whether there exists a point $y_0 \in X(\overline{K})$ satisfying~$\zeta_{\mathrm{abs}}(\overline{D})  = \widehat{h}_{\overline{D}}(y_0)$. 

In this paper, we show that $Y_0$ cannot be chosen to be zero-dimensional in general. Indeed, we construct an example over $K = \QQ$ with $X = \PP^1$, where $\overline{D}$ is a semipositive adelic $\ZZ$-Cartier divisor. The proofs of Theorem 1.2, Theorem 5.1, and Lemma 5.2 in \cite{Ba23} do not seem to yield this information on $Y_0$. Since there is a natural correspondence between semipositive adelic $\ZZ$-Cartier divisors and pairs of a semipositive adelic line bundle and a choice of trivialization (see \cite[Remarks~3.2~and~3.5]{Ba23}), we formulate our theorem in terms of adelic line bundles. 

\begin{theorem}
\label{thm:main}
There exists a semipositive adelic line bundle $\overline{\OO_{\PP^1}(1)} = \left(\OO_{\PP^1}(1), \{\Vert\ndot\Vert_p\}_{p \in M_\QQ}\right)$ on~$\PP^1$ over $\QQ$ with the following properties\textup{:}
\begin{enumerate}
\item[\textup{(i)}]
$\zeta_\mathrm{abs}(\overline{\OO_{\PP^1}(1)}) = 0$.
\item[\textup{(ii)}]
For any $x \in \PP^1(\overline{\QQ})$, one has $\widehat{h}_{\overline{\OO_{\PP^1}(1)}}(x) > 0$. 
\end{enumerate}
\end{theorem}

We note that, by Zhang's fundamental inequality on successive minima (\cite[Theorem~(5.2)]{Zh95a}), 
this example satisfies  
\[
0 = \zeta_\mathrm{abs}(\overline{\OO_{\PP^1}(1)})
= \widehat{h}_{\overline{\OO_{\PP^1}(1)}}(\PP^1) 
= \sup_{
\substack{
F \subset \PP^1(\overline{\QQ})\\  
\text{$F$ finite}}} \inf_{\substack{x \in \PP^1(\overline{\QQ})\\x\notin{F}}} 
\widehat{h}_{\overline{\OO_{\PP^1}(1)}}(x). 
\]
As a result, arithmetic equidistribution (\textit{e.g.}, \cite{Yu08}) is applicable to this setting, even though the set $\{x \in \PP^1(\overline{\QQ}) \mid \widehat{h}_{\overline{\OO_{\PP^1}(1)}}(x) = 0\}$ is empty. 

To construct the semipositive adelic line bundle for Theorem~\ref{thm:main}, we use ideas from arithmetic dynamics. 
Following \cite{Ka07}, we consider a sequence of suitable morphisms $f_n\colon \PP^1 \to \PP^1$ for $n \in \ZZ_{>0}$ and the canonical height function $\widehat{h}_{\mathcal{F}}$ on $\PP^1(\overline{\QQ})$ associated to the sequence $\mathcal{F} \colonequals (f_n)_{n=1}^\infty$. In this situation, if the morphisms $f_n$ are polynomials, then $\widehat{h}_{\mathcal{F}}$ takes the minimum at the point at infinity. We thus need to work with non-polynomial maps. Inspired by \cite{PTW22}, given positive integers $c_n$ and $d_n$, we consider the morphism given by $f_n(X) = \frac{1}{X^{d_n}}+c_n$, which is obtained by the composition of the involution $X \mapsto \frac{1}{X}$ followed by  the polynomial map $X \mapsto X^{d_n}+c_n$. By carefully choosing $(c_n)_{n=1}^\infty$ and $(d_n)_{n=1}^\infty$, we show that 
$\widehat{h}_{\mathcal{F}}$ is the height function associated to a semipositive adelic line bundle on $\PP^1$ that satisfies Properties~(i) and (ii) of Theorem~\ref{thm:main}. 

\smallskip
{\sl Acknowledgments.}\quad
SK thanks Liang-Chung Hsia and Thomas Tucker for helpful discussions during his stay at Academia Sinica in Taiwan, and Julie Tzu-Yueh Wang and the institute for their kind hospitality. We thank Fran\c{c}ois Balla\"{y} for helpful comments. 
SK is supported by KAKENHI 23K03041 and 25H00587. FP is supported by IRN MaDeF and IRN GandA (CNRS).  

\section{Proof of Theorem~\ref{thm:main}}

\subsection{Construction of morphisms}
We embed the affine line $\mathbb{A}^1$ into the projective line $\PP^1$ via the morphism $\alpha \mapsto (\alpha:1)$, and define the point at infinity as $\infty \colonequals (1:0)$. Under this identification, we regard the set of algebraic points as $\PP^1(\overline{\QQ}) =  \overline{\QQ}\cup\{\infty\}$. 

Fix a constant $C > 1$. We will inductively define a sequence of morphisms  
$f_n\colon \PP^1 \to \PP^1$ over $\QQ$. Choosing positive integers $c_n$ and $d_n$ with $d_n \geqslant 2$, we define the rational function 
\[
f_n(X) = \frac{1}{X^{d_n}} + c_n \in \QQ(X), 
\]
which induces the morphism 
\begin{equation}
\label{eqn:def:f:n}
    f_n : \begin{array}{ccc}
            \PP^1 & \longrightarrow & \PP^1 \\
            (x_0:x_1) & \longmapsto & (c_n x_0^{d_n} + x_1^{d_n}: x_0^{d_n}).
          \end{array}
\end{equation}

To begin the induction, 
for $n =1$, we choose any positive integers $c_1$ and $d_1$ satisfying  
$2 \leqslant d_1$ and $\frac{\log (c_1+1)}{d_1} \leqslant C$, and 
set $f_1(X) = \frac{1}{X^{d_1}} + c_1$. Now assume that 
we have defined $f_1(X), \ldots, f_n(X)$ by choosing 
sequences $(c_i)_{i=1}^n$ and $(d_i)_{i=1}^n$ such that  
$c_1 < \cdots < c_n$ and $2 \leqslant d_1 < \cdots < d_n$. 

We construct $f_{n+1}(X)$ by choosing $c_{n+1}$ and $d_{n+1}$ that satisfy the conditions  described in Lemma~\ref{lem:construction} below. 

For $1 \leqslant i \leqslant n$ and $0 \leqslant \ell \leqslant n-i$, 
we define 
\[
S(i, \ell) \colonequals \{\beta \in \overline{\QQ}\cup\{\infty\} \mid (f_{i + \ell} \circ \cdots \circ f_i)(\beta) = \beta\}.
\] 
Note that $\# S(i, \ell) \leqslant \deg(f_{i + \ell} \circ \cdots \circ f_i) +1 \leqslant 
\deg(f_{n} \circ \cdots \circ f_1) + 1 \leqslant d_n \cdots  d_1 +1$.

\begin{lemma}
\label{lem:construction}
There exist positive integers $c_{n+1}$ and $d_{n+1}$ with $c_n < c_{n+1}$, with $d_n < d_{n+1}$, and with $\frac{\log (c_{n+1}+1)}{d_{n+1}} \leqslant C$, such that if we set $f_{n+1}(X) = \frac{1}{X^{d_{n+1}}} + c_{n+1}$, then 
\[
(f_{n+1} \circ \cdots \circ f_{n+1 - \ell})(\beta) \neq \beta
\]
for any $1 \leqslant i \leqslant n$, $0 \leqslant \ell \leqslant n-i$, and 
$\beta \in S(i, \ell)$. 
\end{lemma}

\Proof
We set $S \colonequals 
\bigcup_{1 \leqslant i \leqslant n} \bigcup_{0 \leqslant \ell \leqslant n-i} S(i, \ell)$. 
Then the cardinality of $S$ is bounded above by  
$n \cdot (n+1) \cdot (d_n \cdots  d_1 +1)$. 
For any $d \geqslant 1$, we define a subset $T_d$ of $\overline{\QQ}\cup\{\infty\}$ by 
\[
T_d \colonequals 
\left\{\left.
\beta - \frac{1}{\big((f_{n} \circ \cdots \circ f_{n+1 - \ell})(\beta)\big)^d} 
\;\right|\; \beta \in S
\right\}. 
\]
Here, if $\ell = 0$, then $(f_{n} \circ \cdots \circ f_{n+1 - \ell})(\beta)$ means $\beta$. 
Thus if  $\ell = 0$ and $\beta = \infty$, then $\beta -\big((f_{n} \circ \cdots \circ f_{n+1 - \ell})(\beta)\big)^{-d}$ is equal to $\infty$. If $\ell \geqslant 1$ and $\beta = \infty$, then 
observing that $f_{n+1 - \ell}(\beta) = c_{n+1 - \ell} > 0$, we see that  
$(f_{n} \circ \cdots \circ f_{n+1 - \ell})(\beta)$ is a positive real number. Thus 
if $\ell \geqslant 1$ and $\beta = \infty$, then $\beta - \big((f_{n} \circ \cdots \circ f_{n+1 - \ell})(\beta)\big)^{-d}$ is equal to $\infty$.

Note that $\# T_d \leqslant n \cdot (n+1) \cdot (d_n \cdots  d_1 +1)$. 
We choose a sufficiently large $d$ such that $d_n < d$ and 
$c_n + n \cdot (n+1) \cdot (d_n \cdots  d_1 +1) + 2 < \exp(C d)$. Then 
by Dirichlet's box principle, 
there exists a positive integer $c$ such that $c_n < c$, $c+1 \leqslant \exp(C d)$, and 
$c \not\in T_d$. We set $c_{n+1} \colonequals c$ and $d_{n+1} \colonequals d$, 
and $f_{n+1}(X) = \frac{1}{X^{d_{n+1}}} + c_{n+1}$. Then $f_{n+1}(X)$ satisfies 
the desired properties.
\QED

By induction, we have defined morphisms 
$f_n\colon \PP^1 \to \PP^1$ over $\QQ$ for all $n \in \ZZ_{>0}$. 

\begin{rem}
\label{rmk:lem:construction}
Let $1 \leqslant i_1 < i_2$ and $0 \leqslant \ell$. 
We set $n \colonequals i_2 + \ell - 1$, which implies that $n \geqslant i_1 + \ell$. 
Suppose that 
$\beta \in \overline{\QQ}\cup\{\infty\}$ satisfies 
$\beta = f_{i_1 + \ell} \circ \cdots \circ f_{i_1}(\beta)$. Then 
$\beta \in S(i_1, \ell)$. By Lemma~\ref{lem:construction}, 
we have $\beta \neq f_{i_2 + \ell} \circ \cdots \circ f_{i_2}(\beta)$. 
\end{rem}

\subsection{Height functions}
Let $K$ be a number field and denote by $M_K$ the set of 
places of~$K$. We write $M_K^\infty$ for the set of archimedean places 
and $M_K^0$ for the set of nonarchimedean places. 
For each $v \in M_K$, let $K_v$ be the completion of $K$ with respect to $v$, 
and let $|\ndot|_v$ denote the unique absolute value on $K_v$ extending the usual 
absolute value $|\ndot|_p$ on $\QQ_p$ for $p \mid v$. Specifically, if $p = \infty$, then $|\ndot|_p$ is the usual absolute value on $\RR$; if $p$ is a prime, then $|\ndot|_p$ 
is the $p$-adic absolute value such that $|p|_p = 1/p$. 

The {\em absolute logarithmic Weil height function} $h\colon \PP^1(\overline{\QQ}) \to \RR_{\geqslant 0}$ is defined by 
\[
h(x) \colonequals \sum_{v \in M_K} \frac{[K_v:\QQ_p]}{[K:\QQ]} \log \max \{|x_0|_v, |x_1|_v\}, 
\]
where $x = (x_0:x_1)$ and $K$ is any number field containing $x_0$ and $x_1$. 

\begin{lemma}
\label{lem:eg}
Let $c$ and $d$ be positive integers, and consider the rational function given by
$f(X) =   \frac{1}{X^d} + c\in \QQ(X)$. As in the previous subsection, we regard $f$ as a morphism $f\colon \PP^1 \to \PP^1$ over $\QQ$ defined by 
\begin{equation}
\label{eqn:def:f}
f : \begin{array}{ccc}
            \PP^1 & \longrightarrow & \PP^1 \\
            (x_0:x_1) & \longmapsto & \big(c\, x_0^d + x_1^d: x_0^d\big).
          \end{array}
\end{equation}
Then we have 
\[
\left\Vert h - \frac{1}{d} h \circ f  
\right\Vert_{\sup}
\colonequals 
\sup_{x \in \PP^1(\overline{\QQ})}
\left\vert h(x) - \frac{1}{d} h(f(x)) 
\right\vert
\leqslant \frac{\log (c+1)}{d}. 
\]
\end{lemma}

\Proof
Let $x = (x_0: x_1) \in \PP^1(\overline{\QQ})$. Let $K$ be a number field containing $x_0$ and $x_1$. Fix a place $v \in M_K$. 

Assume first that $v$ is non-archimedean. 
Since $|c|_v \leqslant 1$, the strong triangle inequality gives 
\begin{equation}
\label{eqn:padic}
\log \max\{|c\, x_0^d + x_1^d|_v, \; |x_0^d|_v\}
= \log \max\{|x_1^d|_v, \; |x_0^d|_v\}
= d \log \max \{|x_0|_v, \; |x_1|_v\}. 
\end{equation}

Assume next that $v$ is archimedean. Noting that 
$|c|_v = c$, the usual triangle inequality gives 
\begin{align*}
\log \max\{|c\, x_0^d + x_1^d|_v, |x_0^d|_v\}
& \leqslant \log \max\{c\, |x_0^d|_v + |x_1^d|_v, |x_0^d|_v\} \\
& \leqslant d \log \max \{|x_0|_v, \; |x_1|_v\} + \log (c+1),  \\
d \log \max\{|x_0|_v, \; |x_1|_v\}
& \leqslant \log \max\{|c\, x_0^d + x_1^d|_v + c |x_0^d|_v, |x_0^d|_v\} \\
& \leqslant \log \max\{|c\, x_0^d + x_1^d|_v, |x_0^d|_v\} + \log (c+1). 
\end{align*}
Setting $n_v \colonequals[K_v:\QQ_p]/[K:\QQ]$ and summing over all $v \in M_K$, we obtain 
\begin{multline*}
\left|
\sum_{v \in M_K} n_v
\log \max\{|x_0|_v, |x_1|_v\} 
- \frac{1}{d} \sum_{v \in M_K} n_v
\log \max\{|c\, x_0^d + x_1^d|_v, |x_0^d|_v\} 
\right| 
\leqslant \frac{\log (c+1)}{d}, 
\end{multline*}
as desired. 
\QED

Recall that $C > 1$ is a fixed real number, and 
we have constructed strictly increasing sequences~$(c_n)_{n = 1}^\infty$ and $(d_n)_{n = 1}^\infty$ of positive integers 
satisfying $\log (c_n+1)/d_n \leqslant C$, which further enjoy  
the properties described in Lemma~\ref{lem:construction}. For each~$n \in \ZZ_{>0}$, 
let $f_n\colon \PP^1 \to  \PP^1$ be 
the morphism defined in~\eqref{eqn:def:f:n}. 

We set $\mathcal{F} \colonequals (f_n)_{n=1}^\infty$ and 
define a non-negative function 
$\widehat{h}_{\mathcal{F}}\colon  \PP^1(\overline{\QQ}) \to \RR_{\geqslant 0}$ by 
\begin{equation}
\label{eqn:h:Fcal}
\widehat{h}_{\mathcal{F}}(x) 
\colonequals \lim_{n\to\infty} 
\frac{1}{d_n \cdots d_1} 
h\left((f_n \circ \cdots \circ f_1)(x)\right)
\qquad \text{(for $x \in \PP^1(\overline{\QQ})$)}. 
\end{equation}

Letting $h_n(x) \colonequals 
\frac{1}{d_n \cdots d_1} 
h\left((f_n \circ \cdots \circ f_1)(x)\right)$, we have the following proposition. 

\begin{proposition}[{\cite{CS93}, \cite{Ka07}}]
\label{prop:Ka}
The following properties hold:
\begin{enumerate}
\item
The limit defining $\widehat{h}_{\mathcal{F}}(x) $ in \eqref{eqn:h:Fcal} exists for 
any $x \in \PP^1(\overline{\QQ})$. 
\item
We have $\Vert \widehat{h}_{\mathcal{F}} - h \Vert_{\sup} \leqslant 2C$. 
\item
We have $\Vert \widehat{h}_{\mathcal{F}} - h_n \Vert_{\sup} \leqslant \frac{C}{2^{n-1}}$ for any fixed $n\geqslant1$. 
\item
Let $\alpha \in \PP^1(\overline{\QQ})$. Then 
$\widehat{h}_{\mathcal{F}}(\alpha) = 0$ if and only if the forward orbit 
$\{(f_{n} \circ\cdots\circ f_1)(\alpha) \mid n \geqslant 0\}$ is finite. 
\end{enumerate}
\end{proposition}

\Proof
This statement is a special case of \cite[Theorem~3.3]{Ka07}, which follows from \cite{CS93}. For  the reader's convenience, we sketch a proof. By Lemma~\ref{lem:eg}, we have 
$$\left\vert h(x) - \frac{1}{d_i} h( f_i (x))
\right\vert
\leqslant \frac{\log (c_i+1)}{d_i} \leqslant C.$$ 
Thus we get
\[
\left\vert \frac{1}{d_{n-1} \cdots d_1} 
h\left((f_{n-1} \circ \cdots \circ f_1)(x)\right) - \frac{1}{d_n \cdots d_1} h\left((f_n \circ \cdots \circ f_1)(x)\right)  
\right\vert
\leqslant \frac{1}{d_{n-1} \cdots d_1} C. 
\]
Since we have $\left(1 + \sum_{n=1}^\infty \frac{1}{d_n \cdots d_1}\right) C
\leqslant \left(1 + \sum_{n=1}^\infty \frac{1}{2^n}\right) C 
= 2C$, we see that the sequence
$$\left\{\frac{1}{d_n \cdots d_1} h\left((f_n \circ \cdots \circ f_1)(x)\right)\right\}_{n=1}^\infty$$ is a Cauchy sequence, and 
we obtain (1) and~(2). Furthermore, since 
\[
\frac{1}{d_n \cdots d_1} \left(1 + \sum_{k=n+1}^\infty \frac{1}{d_k \cdots d_{n+1}}\right) C 
\leqslant \frac{1}{2^n} \cdot 2C = \frac{1}{2^{n-1}} C, 
\]
we obtain (3). 

For (4), the ``if'' direction is an immediate consequence of  the definition of $\widehat{h}_{\mathcal{F}}$. To show the ``only if'' direction, we define the shifted sequence $\mathcal{F}_{k}  \colonequals (f_n)_{n=k}^\infty$ for $k \geqslant 1$. For $n \geqslant k$, we set 
\begin{equation}
\label{eqn:h:Fcal:k}
\widehat{h}_{\mathcal{F}_k}(x) 
\colonequals \lim_{n\to\infty} 
\frac{1}{d_n \cdots d_k} 
h\left((f_n \circ \cdots \circ f_k)(x)\right).   
\end{equation}
Applying the same proof used for (2), we obtain  
$\Vert \widehat{h}_{\mathcal{F}_k} - h\Vert_{\sup} \leqslant 2C$. 
Furthermore, comparing 
\eqref{eqn:h:Fcal} and \eqref{eqn:h:Fcal:k} gives 
$\widehat{h}_{\mathcal{F}_k}\left(
(f_{k-1}\circ \cdots\circ f_1)(\alpha) 
\right) 
= d_{k-1}\cdots d_1\, \widehat{h}_{\mathcal{F}}(\alpha) = 0
$. Thus 
\[
\{(f_{n} \circ\cdots\circ f_1)(\alpha) \mid n \geqslant 0\}
\subseteq \{x \in \PP^1(\QQ(\alpha)) \mid 
h(x) \leqslant 2C\}. 
\]
Since the latter set is finite by the Northcott property, the former set is also finite. 
\QED

We call $\widehat{h}_{\mathcal{F}}\colon  \PP^1(\overline{\QQ}) \to \RR_{\geqslant 0}$ the {\em canonical height function} associated to 
$\mathcal{F} = (f_n)_{n=1}^\infty$.

\subsection{Green functions}
In this subsection, 
we consider Green functions related to $(f_n)_{n=1}^\infty$, which will be used for adelic metrics in the next subsection. 

For $n \in \ZZ_{> 0}$, let $f_n\colon \PP^1 \to \PP^1$ be the morphism 
defined in \eqref{eqn:def:f:n}. By base change to $\CC$, $f_n$ induces a morphism $f_n\colon \PP^1_\CC \to \PP^1_\CC$, which 
we still denote by $f_n$. 
Let $Z_0, Z_1$ denote the coordinate functions of $\CC^2 = \mathbb{A}_\CC^2(\CC)$. 
We consider the lift $F_n$ of $f_n$ defined by 

\begin{equation}
\label{eqn:def:f}
F_n : \begin{array}{ccc}
            \mathbb{A}_\CC^2\setminus\{O\} & \longrightarrow & \mathbb{A}_\CC^2\setminus\{O\} \\
            (z_0, z_1) & \longmapsto & (c_n z_0^{d_n} + z_1^{d_n}, z_0^{d_n}).
          \end{array}
\end{equation}
where $O = (0, 0)$ denotes the origin of $\mathbb{A}_\CC^2$. 
We set $\Phi_n = F_n \circ \cdots \circ F_1$, and write 
$\Phi_n = (\Phi_{n, 0}, \Phi_{n, 1})$, where 
$\Phi_{n, 0}, \Phi_{n, 1} \in \CC[Z_0, Z_1]$. 
We then define a function $G_n\colon \mathbb{A}_\CC^2(\CC)\setminus\{O\} \to \RR$ by 
\begin{align*}
G_n(z) 
& = \frac{1}{d_n \cdots d_1} \log 
\sqrt{\left| \Phi_{n, 0}(z) \right|^2 + \left| \Phi_{n, 1}(z) \right|^2}
\\
& = \frac{1}{d_n \cdots d_1} \log 
\big(\left| \Phi_{n, 0}(z) \right|^2 + \left| \Phi_{n, 1}(z) \right|^2\big)^{1/2}
\qquad (z \in \CC^2\setminus\{O\}). 
\end{align*}

For $z = (z_0, z_1)\in  \CC^2\setminus\{O\}$, we set 
\begin{equation}
\label{eqn:def:G:F}
\widehat{G}_{\mathcal{F}}(z) = \lim_{n\to \infty} G_n(z). 
\end{equation}

\begin{proposition}
\label{prop:Green:n}
The following properties hold:
\begin{enumerate}
\item
The limit defining $\widehat{G}_{\mathcal{F}}(z)$ in 
\eqref{eqn:def:G:F} exists for any $z\in  \CC^2\setminus\{O\}$. 
\item
The sequence of functions $G_n$ converges uniformly 
to $\widehat{G}_{\mathcal{F}}$ on $ \CC^2\setminus\{O\}$. 
\end{enumerate}
\end{proposition}

\Proof
We write $F_n= (F_{n, 0}, F_{n, 1})$, where $F_{n, 0}, F_{n, 1} \in \CC[Z_0, Z_1]$. 
We first show that for any $z = (z_0, z_1) \in \CC^2\setminus\{O\}$,   
\begin{equation}
\label{eqn:diff:Fn}
\left|
\log \big(|z_0|^2 + |z_1|^2\big)^{1/2} - 
\frac{1}{d_n}
\log \big(\left|F_{n, 0}(z)\right|^2 + \left|F_{n, 1}(z)\right|^2\big)^{1/2} 
\right|
\leqslant \frac{\log(c_n+1)}{d_n}  + \frac{\log 2}{2}. 
\end{equation}
Indeed, we compute
\begin{align*}
& \log \left(\left|F_{n, 0}(z)\right|^2 + \left|F_{n, 1}(z)\right|^2\right)^{1/2}
= \log \left(|c_n z_0^{d_n} + z_1^{d_n}|^2 + |z_0^{d_n}|^2\right)^{1/2} \\
& \qquad \leqslant 
\log \left(\left(|c_n z_0^{d_n}| + |z_1^{d_n}|\right)^2 + |z_0^{d_n}|^2\right)^{1/2} \leqslant 
 \log (|z_0^{d_n}|^2 + |z_1^{d_n}|^2)^{1/2} + \log(c_n+1) 
\\
& \qquad \leqslant 
d_n \log (|z_0|^2 + |z_1|^2)^{1/2} + \log(c_n+1). 
\end{align*}
Similarly, we compute
\begin{align*}
& d_n \log (|z_0|^2 + |z_1|^2)^{1/2}
\leqslant \log \left(|z_0^{d_n}|^2 +  |z_1^{d_n}|^2\right)^{1/2} + \frac{d_n}{2} \log 2 \\
& \qquad \leqslant 
\log  \left(|z_0^{d_n}|^2 + \left(|c_n z_0^{d_n} + z_1^{d_n}| +  |c_n z_0^{d_n}|\right)^2\right)^{1/2} +  \frac{d_n}{2} \log 2 \\
& \qquad  \leqslant \log \left(|c_n z_0^{d_n} + z_1^{d_n}|^2 + |z_0^{d_n}|^2\right)^{1/2}
+ \log (c_n+1) +  \frac{d_n}{2} \log 2. 
\end{align*}
Dividing by $d_n$ gives \eqref{eqn:diff:Fn}. 

Once we confirm \eqref{eqn:diff:Fn}, the remainder of the proof  follows \cite[\S2]{FW00} or  \cite[Proposition~6.3]{Ka07}. Recall from the construction of $f_n$ that $C > 1$ is a constant satisfying $\frac{\log(c_n+1)}{d_n} \leqslant C$ for all $n$. We set $C^\prime \colonequals C + (\log 2)/2$. Replacing $z$ with $\Phi_{n-1}(z)$ in \eqref{eqn:diff:Fn}  and dividing by $d_{n-1}\cdots d_1$, we find 
\[
\left|
G_{n-1}(z) 
- G_{n}(z) 
\right| 
\leqslant \frac{1}{d_{n-1}\cdots d_1}C^\prime
\]
for any $z = (z_0, z_1) \in \CC^2\setminus\{O\}$. 

By the same argument used in the proof of Proposition~\ref{prop:Ka}, 
$\{G_n(z)\}_{n=1}^\infty$ is a Cauchy sequence, and we have the bound 
$|G_n(z) - \widehat{G}_{\mathcal{F}}(z)| \leqslant C^\prime/2^{n-1}$. 
We obtain (1) and~(2).  
\QED

We call $\widehat{G}_{\mathcal{F}}\colon \mathbb{A}^2(\CC)\setminus \{O\}\to \RR$ the {\em Green function} associated to $\mathcal{F} = (f_n)_{n=1}^\infty$. 

Noting that $G_{n}(\lambda z) = G_{n}(z) + \log |\lambda|$ and $\widehat{G}_{\mathcal{F}}(\lambda z) = \widehat{G}_{\mathcal{F}}(z)+ \log |\lambda|$ for $\lambda \in \CC\setminus\{0\}$, 
we define the functions $g_{n}$ and $\widehat{g}_{\mathcal{F}}$ on 
$\PP^1(\CC)$ by 
\begin{align}
\label{eqn:def:g:F}
g_{n}(x)
\colonequals G_{n}(z)  - \log \sqrt{|z_0|^2 + |z_1|^2}, 
\\
\widehat{g}_{\mathcal{F}}(x)
\colonequals \widehat{G}_{\mathcal{F}}(z)  - \log \sqrt{|z_0|^2 + |z_1|^2}, 
\end{align}
where $x = (x_0: x_1) \in \PP^1(\CC)$ and $z = (x_0, x_1) \in \CC^2\setminus\{O\}$ 
is any lift of the point $x$. These functions are well-defined and independent of the choice of lift $z$. 

We define metrics $\Vert\ndot\Vert_n$ and $\Vert\ndot\Vert^{\mathcal{F}}$ 
on $\OO_{\PP^1}(1)$ by 
$\Vert\ndot\Vert_n = \exp(-g_n)\, 
\Vert\ndot\Vert_{FS}$ and 
$\Vert\ndot\Vert^{\mathcal{F}} = \exp(-\widehat{g}_{\mathcal{F}})\, 
\Vert\ndot\Vert_{FS}$, where 
$\Vert\ndot\Vert_{FS}$ is the Fubini-Study metric: 
For $s = a_0 X_0 + a_1 X_1 \in H^0\big(\PP^1, \OO_{\PP^1}(1)\big)$, 
\[
\Vert s\Vert_{FS}(x) = \frac{|a_0 x_0 + a_1 x_1|}{\sqrt{|x_0|^2 + |x_1|^2}}
\qquad (x = (x_0:x_1) \in \PP^1(\CC)). 
\]

\begin{lemma}
\label{lem:FS}
The following properties hold:
\begin{enumerate}
\item
The metrics $\Vert\ndot\Vert_n$ converge uniformly 
to $\Vert\ndot\Vert^{\mathcal{F}}$ on $\PP^1(\CC)$ as $n \to \infty$. 
\item
For each $n$, $\Vert\ndot\Vert_n$ is smooth, and its first Chern form 
$c_1(\OO_{\PP^1}(1), \Vert\ndot\Vert_n)$ is positive. 
\end{enumerate}
\end{lemma}

\Proof
Assertion~(1) follows from Proposition~\ref{prop:Green:n}~(2). 

For (2), since $G_n$ is a $C^\infty$-function on $\CC^2 \setminus \{O\}$, the metric $\Vert\ndot\Vert_n$ is smooth. For a section $s = a_0 X_0 + a_1 X_1 \in  H^0\big(\PP^1, \OO_{\PP^1}(1)\big)$ and a point $x = (x_0:x_1) \in \PP^1(\CC)$, we have $-\log\Vert s \Vert_n(x) = G_n(x_0, x_1) - \log |a_0 x_0 + a_1 x_1|$. Thus  
$c_1(\OO_{\PP^1}(1), \Vert\ndot\Vert_n) = dd^c(2 G_n) = (f_n \circ \cdots \circ f_1)^* \omega_{FS}$, where $\omega_{FS}$ is the Fubini-Study form. Thus $c_1(\OO_{\PP^1}(1), \Vert\ndot\Vert_n)$ is positive, completing the proof of~(2). 
\QED

We also call $\widehat{g}_{\mathcal{F}}$ the {\em Green function} associated to $\mathcal{F}$. 

\subsection{Semipositive adelic metrics}
Let $\widehat{h}_{\mathcal{F}}\colon  \PP^1(\overline{\QQ}) \to \RR_{\geqslant 0}$ be the canonical height function associated to 
$\mathcal{F} = (f_n)_{n=1}^\infty$ defined in \eqref{eqn:h:Fcal}. 

In this subsection, we show that there exists a semipositive adelic line bundle $\overline{\OO_{\PP^1}(1)}^{\mathcal{F}} = \left(\OO_{\PP^1}(1), \{\Vert\ndot\Vert_p^{\mathcal{F}}\}_{p \in M_\QQ}\right)$ on~$\PP^1$ over $\QQ$ such that $\widehat{h}_{\overline{\OO_{\PP^1}(1)}^\mathcal{F}} = \widehat{h}_{\mathcal{F}}$. Here we use the superscript $\mathcal{F}$ to emphasize that 
the metrics depend on $\mathcal{F}$. 

For $p \in M_\QQ$, denote by $\QQ_p$ the completion of $\QQ$ with respect to $|\ndot|_p$ and let $\overline{\QQ}_p$ be a fixed algebraic closure of $\QQ_p$. We also let $|\ndot|_p$ denote the unique absolute value on $\overline{\QQ}_p$ extending $|\ndot|_p$. 
Following \cite{Zh95b} (see also \cite[Appendix~A]{YZ26}), we recall the notion of semipositive adelic line bundles. Here we restrict ourselves to line bundles on $\PP^1$ over $\QQ$, and we let $L$ be such a line bundle. 

We first recall the notion of metrics on line bundles over $\QQ_p$. 
Let  $L_{\QQ_p} \colonequals L \otimes_\QQ \QQ_p$ be the induced line bundle on $\PP^1_{\QQ_p} \colonequals  \PP^1 \times_{\Spec \QQ} {\Spec \QQ_p}$. 
By a {\em $\QQ_p$-metric} $\Vert\ndot\Vert_p$, we mean a $\mathrm{Gal}(\overline{\QQ}_p/\QQ_p)$-invariant collection of $\overline{\QQ}_p$-norms on the fibers $L_{\QQ_p}(x)$ for each $x \in \PP^1_{\QQ_p}(\overline{\QQ}_p)$. 

Assume that $p \in M_\QQ^\infty$, \textit{i.e.}, $p = \infty$. A $\QQ_p$-metric is {\em semipositive} if it is the uniform limit 
of a sequence of smooth metrics on $L_{\QQ_p}$ with semipositive 
curvature forms. 

Assume that $p \in M_\QQ^0$, \textit{i.e.}, $p$ is a prime. 
Suppose that $(\mathcal{X}, \mathcal{M})$ is a {\em projective model} of $(\PP^1_{\QQ_p}, L_{\QQ_p}^{\otimes e})$ over $\ZZ_p$ for some positive integer $e$. That is, $\mathcal{X}$ is a flat and projective integral scheme over $\ZZ_p$ with generic fiber isomorphic to $\PP^1_{\QQ_p}$ and $\mathcal{M}$ is a line bundle on $\mathcal{X}$ whose generic fiber is isomorphic to $L_{\QQ_p}^{\otimes e}$. Any point $x \in  \PP^1_{\QQ_p}(\overline{\QQ}_p)$ extends to a point $\overline{x} \in \mathcal{X}(\overline{\ZZ}_p)$ by taking the Zariski closure. Then a $\overline{\QQ}_p$-norm $\Vert\ndot\Vert$ on $L_{\QQ_p}(x)$ is defined  for $s \in L_{\QQ_p}(x)$ by $\Vert s \Vert = \inf \{|\ell|_p^{1/e} \mid \ell \in \overline{\QQ}_p, \, s^{\otimes e} \in \ell \mathcal{M}(\overline{x})\}$. 
Patching these norms together by varying $x\in \PP^1_{\QQ_p}(\overline{\QQ}_p)$, we obtain a $\QQ_p$-metric of $L_{\QQ_p}$ on $\PP^1_{\QQ_p}(\overline{\QQ}_p)$, which is called the {\em model metric} induced by $(\mathcal{X}, \mathcal{M})$. 
The model metric is {\em semipositive} if $\deg_{\mathcal{M}}(\mathcal{Z}) \geqslant 0$ for any projective curve $\mathcal{Z}$ in the special fiber of~$\mathcal{X}$. 

Keeping the assumption that $p \in M_k^0$, 
a $\QQ_p$-metric $\Vert\ndot\Vert_p$ is {\em semipositive} if it is the uniform limit of a sequence of semipositive model metrics $\Vert\ndot\Vert_{n}$ on $L_{\QQ_p}$, in the sense that the function $\Vert\ndot\Vert_{n}/\Vert\ndot\Vert$ converges uniformly to $1$ on $\PP^1_{\QQ_p}(\overline{\QQ}_p)$. 

\smallskip
Now we define a semipositive adelic line bundle. A {\em semipositive adelic line bundle} $\overline{L} = \left(L, \{\Vert\ndot\Vert_{p}\}_{p \in M_\QQ}\right)$ on $\PP^1$ consists of a line bundle $L$ and a collection $\{\Vert\ndot\Vert_{p}\}_{p \in M_\QQ}$ satisfying the following properties:
\begin{enumerate}
\item[(i)]
For each $p \in M_\QQ$, $\Vert\ndot\Vert_p$ is a semipositive $\QQ_p$-metric 
of $L_{\QQ_p}$ on $\PP^1_{\QQ_p}$. 
\item[(ii)]
There exists a Zariski open subset $U \colonequals  \Spec\ZZ[1/N]$ (for some $N \in \ZZ_{> 0}$) of $\Spec\ZZ$, a projective and flat integral model $\mathcal{X} \to U$ of $X$ over $U$, and a line bundle $\mathcal{L}$ on $\mathcal{X}_U$ extending $L$ on $\PP^1$ such that for each prime $p \in U$, the $\QQ_p$-metric $\Vert\ndot\Vert_p$ is the model metric 
induced by the integral model $(\mathcal{X} \times_{U} \ZZ_{p}, \mathcal{L} \otimes_{\ZZ[1/N]} \ZZ_{p})$.  
\end{enumerate}

Let $\overline{L} = \left(L, \{\Vert\ndot\Vert_{p}\}_{p \in M_\QQ}\right)$ be a semipositive adelic line bundle on $\PP^1$. For $x \in \PP^1(\overline{\QQ})$, take a number field $K$ over which $x$ is defined. For each place $v \in M_K$, taking $p \in M_\QQ$ with 
$v \mid p$, we fix an embedding $\sigma_v\colon K \hookrightarrow \overline{\QQ}_p$. 
Let $s$ be a rational section of $L$ with $s(x) \neq 0$. 
By the $\mathrm{Gal}(\overline{\QQ}_p/\QQ_p)$-invariance of $\Vert\ndot\Vert$, the quantity 
$\Vert s\Vert_v(x) \colonequals \Vert s \Vert_p(\sigma_v(x))$ is independent of the choice of $\sigma_v$. The {\em height} $\widehat{h}_{\overline{L}}(x)$ of $x\in \PP^1(K)$ associated to $\overline{L}$ is defined by 
\begin{equation}
\label{eqn:heiggt:Lbar}
\widehat{h}_{\overline{L}}(x) 
\colonequals 
- \sum_{v \in M_K} 
\frac{[K_v:\QQ_p]}{[K:\QQ]}\,
\log \Vert s \Vert_v(x). 
\end{equation}
By the product formula, 
$\widehat{h}_{\overline{L}}(x)$ is independent of the choice of $s$. 

\smallskip
We now return to the proof of Theorem~\ref{thm:main}. 
For each positive integer $n$, we define 
homogeneous polynomials 
\begin{equation}
\label{eqn:Fn0:Fn1}
F_{n, 0}(X_0, X_1) \colonequals c_n X_0^{d_n} + X_1^{d_n}, \quad 
F_{n, 1}(X_0, X_1) \colonequals X_0^{d_n} \in \ZZ[X_0, X_1]. 
\end{equation}
The morphism 
$f_n\colon \PP^1 \to \PP^1$ over $\QQ$ defined in \eqref{eqn:def:f:n} is given by 
$f_n = (F_{n, 0}: F_{n,1})$. Since the ideal $(F_{n, 0}(X_0, X_1), F_{n, 1}(X_0, X_1)) 
$ is equal to $(X_0^{d_n}, X_1^{d_n})$ in $\ZZ[X_0, X_1]$, $f_n$ extends to a morphism $f_n\colon \PP_\ZZ^1 \to \PP_\ZZ^1$ over $\ZZ$, which we still denote by $f_n$. 

Let  $F_n\colon \mathbb{A}^2_\ZZ \to 
\mathbb{A}^2_\ZZ$ be the morphism given by 
$F_n = (F_{n, 0}, F_{n, 1})$. We set $\Phi_n = F_n \circ \cdots \circ F_1$ and write $\Phi_n = (\Phi_{n, 0}, \Phi_{n, 1})$.  In these coordinates, we have $\phi_n = (\Phi_{n, 0}: \Phi_{n, 1})$. 

Note that $(\PP^1_\ZZ, \OO_{\PP^1_\ZZ}(1))$ is a model 
of $(\PP^1, \OO_{\PP^1}(1))$ over $\Spec \ZZ$ and that 
\begin{equation}
\label{eqn:model}
\left( \PP^1_\ZZ, \; (\phi_n)^*(\OO_{\PP^1_\ZZ}(1))\right)
\end{equation}
 is a model of $(\PP^1, \OO_{\PP^1}(1)^{\otimes D_n})$ over $\Spec \ZZ$. 
We fix the isomorphism 
\begin{equation}
\label{eqn:model:isom}
 \psi_n\colon 
 \OO_{\PP^1_\ZZ}(1)^{\otimes D_n} \overset{\cong}{\longrightarrow}
 (\phi_n)^*(\OO_{\PP^1_\ZZ}(1))
\end{equation}
such that $\psi_n(\Phi_{n, i}) = (\phi_n)^*X_i$ for $i = 0, 1$. 
This defines a model metric on $\OO_{\PP^1_{\QQ_p}}(1)$ 
induced by \eqref{eqn:model} and \eqref{eqn:model:isom} for each $p \in M_\QQ^0$. 

\begin{proposition}
\label{prop:std:metric}
For each $p \in M_\QQ^0$, the model metric 
induced by \eqref{eqn:model}  and \eqref{eqn:model:isom} on $\OO_{\PP^1_{\QQ_p}}(1)$ is the standard metric $\Vert \ndot\Vert_{p, std}$\textup{:} 
for a section $s = a_0 X_0 + a_1 X_1 \in H^0\big(\PP^1_{\QQ_p}, \OO_{\PP^1_{\QQ_p}}(1)\big)$, 
\[
\Vert s\Vert_{p, std}(x) = \frac{|a_0 x_0 + a_1 x_1|_p}{\max\{|x_0|_p, |x_1|_p\}}
\qquad (x = (x_0:x_1) \in \PP^1_{\QQ_p}(\overline{\QQ}_p)). 
\]
\end{proposition}

\Proof
By \cite[(2.3)]{Zh95b}, the induced model metric $\Vert \ndot \Vert$ is the unique metric 
on $\OO_{\PP^1_{\QQ_p}}(1)$ satisfying 
\begin{equation}
\label{eqn:model:metric}
\Vert \ndot \Vert^{D_n} = \psi_n^* \phi_n^* \Vert \ndot \Vert. 
\end{equation}
Let us verify that $\Vert \ndot\Vert_{p, std}$ satisfies \eqref{eqn:model:metric}. Let $s = a_0 X_0 + a_1 X_1 \in H^0\big(\PP^1_{\QQ_p}, \OO_{\PP^1_{\QQ_p}}(1)\big)$ and let 
$x = (x_0:x_1) \in \PP^1(\overline{\QQ}_p)$. We define $t =  \psi_n^* \phi_n^* s \in 
H^0\big(\PP^1_{\QQ_p}, \OO_{\PP^1_{\QQ_p}}(1)^{\otimes D_n}\big)$, which means 
$t = a_0 \Phi_{n, 0}(X_0, X_1) + a_1 \Phi_{n, 1}(X_0, X_1)$.  
We denote the metric on $ \OO_{\PP^1_{\QQ_p}}(1)^{\otimes D_n}$ induced by $\Vert \ndot\Vert_{p, std}$ on the left-hand side of \eqref{eqn:model:metric} by $\Vert\ndot\Vert_{LHS}$  and the right-hand side by $\Vert\ndot\Vert_{RHS}$. 
We have 
$
\Vert s^{D_n} \Vert_{LHS}(x) 
= 
{|(a_0 x_0 + a_1 x_1)^{D_n}|_p}/
{\max\{|x_0|_p^{D_n}, |x_1|_p^{D_n}\}}
$, so 
\begin{align*}
\Vert t \Vert_{LHS}(x) 
= 
\frac{|a_0 \Phi_{n,0}(x_0, x_1) + a_1 \Phi_{n,1}(x_0, x_1)|_p}%
{\max\{|x_0|_p^{D_n}, |x_1|_p^{D_n}\}}. 
\end{align*}
Also, we have 
\begin{align*}
\Vert t \Vert_{RHS}(x) 
= 
\frac{|a_0 \Phi_{n,0}(x_0, x_1) + a_1 \Phi_{n,1}(x_0, x_1)|_p}%
{\max\{|\Phi_{n,0}(x_0, x_1)|_p, |\Phi_{n,1}(x_0, x_1)|_p\}}. 
\end{align*}
Using \eqref{eqn:padic}, we have 
$\max\{|\Phi_{n,0}(x_0, x_1)|_p, |\Phi_{n,1}(x_0, x_1)|_p\}
= \max\{|x_0|_p^{D_n}, |x_1|_p^{D_n}\}$ by induction. 
Thus, we obtain $\Vert\ndot\Vert_{LHS} =\Vert\ndot\Vert_{RHS}$, completing the proof. 
\QED

We define a semipositive adelic line bundle 
\begin{equation}
\label{eqn:O:F}
\overline{\OO_{\PP^1}(1)}^{\mathcal{F}} = \left(\OO_{\PP^1}(1), \{\Vert\ndot\Vert_p^{\mathcal{F}}\}_{p \in M_\QQ}\right)
\end{equation}
on~$\PP^1$ over $\QQ$ as follows: 
\begin{itemize}
\item
For $p = \infty \in M_\QQ^\infty$, we set 
$\Vert\ndot\Vert_\infty^{\mathcal{F}} \colonequals \exp(-\widehat{g}_{\mathcal{F}})\, 
\Vert\ndot\Vert_{FS}$, where $\widehat{g}_{\mathcal{F}}$ is the Green function 
associated to $\mathcal{F}$ in the previous subsection and 
$\Vert\ndot\Vert_{FS}$ is the Fubini-Study metric. 
By Lemma~\ref{lem:FS}, $\Vert\ndot\Vert_\infty^{\mathcal{F}}$ is semipositive. 
\item
For a prime $p$, let $\Vert\ndot\Vert_p^{\mathcal{F}}$ be the standard metric $\Vert \ndot\Vert_{p, std}$. By Proposition~\ref{prop:std:metric},  $\Vert \ndot\Vert_{p, std}$ is semipositive, and thus $\overline{\OO_{\PP^1}(1)}^{\mathcal{F}}$ is indeed an adelic line bundle. 
\end{itemize}

\begin{proposition}
\label{prop:heights:agree}
We have $\widehat{h}_{\mathcal{F}} 
= \widehat{h}_{\overline{\OO_{\PP^1}(1)}^{\mathcal{F}}}$. 
In other words, the height defined in~\eqref{eqn:h:Fcal} coincides 
with the height associated to the semipositive 
adelic line bundle $\overline{\OO_{\PP^1}(1)}^{\mathcal{F}}$. 
\end{proposition}

\Proof
Let $x = (x_0:x_1) \in \PP^1(\overline{\QQ})$, and 
we take a number field $K$ containing $x_0$ and $x_1$. 
Let $s = a_0 X_0 + a_1 X_1 \in H^0(\PP^1, \OO_{\PP^1}(1))$ be a global section such that  
$s(x) \neq 0$. Let $F_n$ be as defined in \eqref{eqn:Fn0:Fn1}. 
Recall that 
$D_n = d_n \cdots d_1$, $\phi_n = f_n \circ \cdots \circ f_1$, and 
$\Phi_n = (\Phi_{n, 0}, \Phi_{n, 1}) \colonequals F_n \circ \cdots \circ F_1$. 
Setting $n_v \colonequals[K_v:\QQ_p]/[K:\QQ]$ and 
using \eqref{eqn:padic}, we compute 
\begin{align*}
& h_n(x) 
\colonequals 
\frac{1}{D_n} 
h\left(\phi_n(x)\right)
\\
&\quad  = 
\sum_{v \in M_K^0} n_v \log \max \{|x_0|_v, |x_1|_v\} 
+ \frac{1}{D_n} \sum_{v \in M_K^\infty}
n_v 
\log \max \{|\Phi_{n, 0}(x_0, x_1)|_v, 
|\Phi_{n, 1}(x_0, x_1)|_v\} \\
& \quad = 
- \sum_{v \in M_K^0} n_v \log \frac{|a_0 x_0 + a_1 x_1|_v}{\max \{|x_0|_v, |x_1|_v\}}
- \frac{1}{D_n} \sum_{v \in M_K^\infty}
n_v \log \frac{|a_0 x_0 + a_1 x_1|_v}{\max \{|\Phi_{n, 0}(x_0, x_1)|_v, 
|\Phi_{n, 1}(x_0, x_1)|_v\}}, 
\end{align*}
where we use the product formula for the last equality. 

For each place $v \in M_K$, we take $p \in M_\QQ$ with $v \mid p$ and fix an embedding $\sigma_v\colon K \hookrightarrow \overline{\QQ}_p$. We set $\Vert\ndot\Vert^{\mathcal{F}}_v(x) \colonequals \Vert \ndot\Vert_{p, std}(\sigma_v(x))$ for $v \in M_K^0$ and 
$\Vert\ndot\Vert^{\mathcal{F}}_v(x) \colonequals \Vert \ndot\Vert_{\infty}^{\mathcal{F}}(\sigma_v(x))$ for $v \in M_K^\infty$. 
Note that for $z = (z_0, z_1) \in \CC^2$, we have 
$\lim_{n\to \infty} \frac{1}{D_n} 
\log \max \{|\Phi_{n, 0}(z_0, z_1)|, |\Phi_{n, 1}(z_0, z_1)|\}
= \lim_{n\to \infty} \frac{1}{D_n} 
\log (|\Phi_{n, 0}(z_0, z_1)|^2 +  |\Phi_{n, 1}(z_0, z_1)|^2)^{1/2}
= \widehat{G}_{\mathcal{F}}(z)$. 
Taking the limit as $n \to \infty$, we obtain 
\[
\widehat{h}_{\mathcal{F}}(x) 
= 
- \sum_{v \in M_K^0} n_v \log \Vert s\Vert^{\mathcal{F}}_v(x)
- \sum_{v \in M_K^\infty}
n_v \log \Vert s\Vert^{\mathcal{F}}_v(x)
= \widehat{h}_{\overline{\OO_{\PP^1}(1)}^{\mathcal{F}}}(x), 
\]
which completes the proof. 
\QED

\subsection{Absolute minima}
Let $\widehat{h}_{\mathcal{F}}\colon  \PP^1(\overline{\QQ}) \to \RR_{\geqslant 0}$ be the canonical height function associated to 
$\mathcal{F} = (f_n)_{n=1}^\infty$ defined in \eqref{eqn:h:Fcal}, 
which by Proposition~\ref{prop:heights:agree} coincides with the height function associated to the semipositive adelic line bundle $\overline{\OO_{\PP^1}(1)}^{\mathcal{F}}$ defined 
in \eqref{eqn:O:F}. 

In this subsection, 
we show that $ \widehat{h}_{\mathcal{F}} = \widehat{h}_{\overline{\OO_{\PP^1}(1)}^\mathcal{F}}$ satisfies Properties (i) and~(ii) in Theorem~\ref{thm:main}. The next proposition gives Property~(i): 
$\zeta_\mathrm{abs}\left(\overline{\OO_{\PP^1}(1)}^\mathcal{F}\right) = 0$. 

\begin{proposition}
We have $\inf_{x \in \PP^1(\QQ)}\widehat{h}_{\mathcal{F}}(x) = 0$. 
\end{proposition}

\Proof
For any $x \in \PP^1(\QQ)$, we have $\widehat{h}_{\mathcal{F}}(x) \geqslant  0$ by \eqref{eqn:h:Fcal}, so $\inf_{x \in \PP^1(\QQ)}\widehat{h}_{\mathcal{F}}(x) \geqslant 0$.

Recall that $h_n(x) \colonequals 
\frac{1}{d_n \cdots d_1} h\left((f_n \circ \cdots \circ f_1)(x)\right)$. 
By Theorem~\ref{prop:Ka}~(3), 
$\Vert \widehat{h}_{\mathcal{F}} -h_n \Vert_{\sup} \to 0$ as~$n \to \infty$. 
Let $\varepsilon > 0$ be any small number. We take a positive  integer $n$ such that $\Vert \widehat{h}_{\mathcal{F}} -h_n \Vert_{\sup} < \varepsilon$. Noting that $f_n \circ \cdots \circ f_1$ is surjective, we take $\alpha_n \in \PP^1(\overline{\QQ})$ such that 
$(f_n \circ \cdots \circ f_1)(\alpha_n) = 0$. We get 
\[
h_n(\alpha_n)
 = \frac{1}{d_n \cdots d_1} h\left((f_n \circ \cdots \circ f_1)(\alpha_n)\right)
= \frac{1}{d_n \cdots d_1} h(0) = 0. 
\]
Then $\widehat{h}_{\mathcal{F}}(\alpha_n) < h_n(\alpha_n) + \varepsilon = \varepsilon$. 
By Proposition~\ref{prop:property:II} below, we have 
$\widehat{h}_{\mathcal{F}}(\alpha_n) > 0$.
It follows that there exist a sequence $(\alpha_n)_{n=1}^\infty$ 
with $\alpha_n \in \overline{\QQ}$ 
and being mutually distinct 
such that $\lim_{n\to\infty} \widehat{h}_{\mathcal{F}}(\alpha_n) = 0$, so $\inf_{x \in \PP^1(\QQ)}\widehat{h}_{\mathcal{F}}(x) = 0$. 
\QED

To show Property~(ii) in Theorem~\ref{thm:main}, we use the following simple lemma. 

\begin{lemma}
\label{lem:elem}
Let $(x_i)_{i = 1}^\infty$ be a sequence with $x_i \in \PP^1(\overline{\QQ})$. Assume that 
$\{x_i \mid i \in \ZZ_{> 0}\}$ is a finite set. 
Then there exist $1 \leqslant i_1 < i_2$ and $0 \leqslant \ell$ such that 
\[
x_{i_1} = x_{i_1 + \ell + 1} = x_{i_2} = x_{i_2 + \ell+1}. 
\]
\end{lemma}

\Proof
We set $\{x_i \mid i \in \ZZ_{> 0}\} = \{y_1, \ldots, y_m\}$. 
Fix any $i \in \ZZ_{> 0}$.  By the Dirichlet box principle, 
there exist integers $k_1, k_2$ with $0 \leqslant k_1 < k_2 \leqslant m$ such that 
$x_{i+k_1} = x_{i+k_2}$. 
It follows that there exists $y  \in \overline{\QQ}$ (which is one 
of $y_1, \ldots, y_m$) such that 
there are increasing sequences $(i_j)_{j=1}^{\infty}$ and $(i_j^\prime)_{j=1}^{\infty}$of positive integers such that 
$1 \leqslant  i_j^\prime - i_j \leqslant m$ and $x_{i_j} = x_{i_j^\prime} = y$. 

Replacing $(i_j)_{j=1}^{\infty}$ by a subsequence, there exists a non-negative integer $0 \leqslant \ell \; (\leqslant m-1)$ such that 
$x_{i_j} = x_{i_j +\ell + 1} = y$ for all $j$. 
Considering the first two indices $i_1$ and $i_2$, 
we have $x_{i_1} = x_{i_1 + \ell+1} = x_{i_2} = x_{i_2 + \ell+1}$. 
\QED

We show Property~(ii) in Theorem~\ref{thm:main}. 

\begin{proposition}
\label{prop:property:II}
For any $\alpha \in \PP^1(\overline{\QQ})$, we have 
$\widehat{h}_{\mathcal{F}}(\alpha) > 0$. 
\end{proposition}

\Proof
To derive a contradiction, 
we assume that there exists $\alpha \in \PP^1(\overline{\QQ})$ with $\widehat{h}_{\mathcal{F}}(\alpha) = 0$. 
By Theorem~\ref{prop:Ka}~(4), $\{(f_{i-1}\circ \cdots\circ f_1)(\alpha) \mid 
i \in \ZZ_{>0}\}$ is a finite set. 

Using Lemma~\ref{lem:elem}, there exist $1 \leqslant i_1 < i_2$ and $0 \leqslant \ell$ such that 
\begin{align*}
(f_{i_1-1}\circ \cdots\circ f_1)(\alpha) 
& = (f_{i_1+ \ell}\circ \cdots\circ f_1)(\alpha) \\
& = (f_{i_2-1}\circ \cdots\circ f_1)(\alpha) 
=  (f_{i_2+ \ell}\circ \cdots\circ f_1)(\alpha). 
\end{align*}
We set $\beta \colonequals (f_{i_1-1}\circ \cdots\circ f_1)(\alpha)
= (f_{i_2-1}\circ \cdots\circ f_1)(\alpha)$. 
Then we have 
\[
\beta = (f_{i_1+ \ell}\circ \cdots\circ f_{i_1})(\beta) 
=  (f_{i_2+ \ell}\circ \cdots\circ f_{i_2})(\beta). 
\]
However, this contradicts the properties of $f_n$ described in Lemma~\ref{lem:construction} (see Remark~\ref{rmk:lem:construction}). This completes the proof. 
\QED

\end{document}